\documentclass[a4paper,11pt]{article}

\usepackage{mathrsfs}
\usepackage{amsmath}
\usepackage[applemac]{inputenc}
\usepackage{amsfonts}
\usepackage{amssymb}
\usepackage{amsthm}
\usepackage{stmaryrd}
\usepackage{color}
\usepackage{mathabx}

\usepackage{cancel}                                                             
\usepackage[normalem]{ulem}

\newcommand{\N}{\mathbb{N}}

\newcommand{\R}{\mathbb{R}}

\def\1{\mbox{1\hspace{-0.25em}l}}

\def \t {\tau}

\def\bx{{\mathbf{x}}}
\def\x{{\boldsymbol{x}}}
\def\p{{\mathbf{p}}}

\def\l{\lambda}
\def\<{\langle}
\def\>{\rangle}

\def\x{{\xi}}
\def\z{{\zeta}}

\numberwithin{equation}{section}

\newcommand{\be}{\begin{eqnarray}}
\newcommand{\ee}{\end{eqnarray}}
\newcommand{\ce}{\begin{eqnarray*}}
\newcommand{\de}{\end{eqnarray*}}
\newtheorem{theorem}{Theorem}[section]
\newtheorem{lemma}[theorem]{Lemma}
\newtheorem{remark}[theorem]{Remark}
\newtheorem{definition}[theorem]{Definition}
\newtheorem{proposition}[theorem]{Proposition}
\newtheorem{assumption}[theorem]{Assumption}
\newtheorem{Examples}[theorem]{Example}
\newtheorem{corollary}[theorem]{Corollary}

\def\eps{\varepsilon}
\def\t{\tau}
\def\e{\mathrm{e}}

\def\a{\alpha}

\def\p{\partial}
\def\d{\delta}
\def\g{\gamma}
\def\l{\lambda}

\def\[{{\Big[}}
\def\]{{\Big]}}
\def\<{{\langle}}
\def\>{{\rangle}}
\def\({{\big(}}
\def\){{\big)}}

\def\bx{{\mathbf{x}}}

\def\sgn{\mbox{\rm sgn}}

\def\bb2{{\boldsymbol{2}}}

\def\={&\!\!=\!\!&}

\def\cG{{\mathcal G}}

\def\1{{\mathbf{1}}}

\def\geq{\geqslant}
\def\leq{\leqslant}
\def\ge{\geqslant}
\def\le{\leqslant}

\def\eps{\varepsilon}
\def\t{\tau}
\def\e{\mathrm{e}}

\def\a{\alpha}

\def\p{\partial}
\def\g{\gamma}
\def\l{\lambda}

\def\[{{\Big[}}
\def\]{{\Big]}}
\def\<{{\langle}}
\def\>{{\rangle}}

\def\bx{{\mathbf{x}}}

\def\sgn{\mbox{\rm sgn}}

\def\={&\!\!=\!\!&}
\def\bt{\begin{theorem}}
\def\et{\end{theorem}}
\def\bl{\begin{lemma}}
\def\el{\end{lemma}}
\def\br{\begin{remark}}
\def\er{\end{remark}}
\def\bx{\begin{Examples}}
\def\ex{\end{Examples}}
\def\bd{\begin{definition}}
\def\ed{\end{definition}}
\def\bp{\begin{proposition}}
\def\ep{\end{proposition}}
\def\bc{\begin{corollary}}
\def\ec{\end{corollary}}

\def\geq{\geqslant}
\def\leq{\leqslant}
\def\ge{\geqslant}
\def\le{\leqslant}

 \def\R{\mathbb R}
 \def\R{\mathbb R}    
\def\N{\mathbb N}  
   
\def\<{\langle} \def\>{\rangle}

\def\hD{\hat{D}}
\def\K{{\mathcal{K}}}

\def\0{{\mathbf{0}}}

\def\G{\Gamma}

\def\b{\beta}
\def\bbB{\langle\b\rangle_B}
\newcommand{\Q}{\mathbf{d}}
\newcommand{\norm}[1]{\left\|{#1}\right\|}
\def\La{\Lambda}
\def\DB{{\mathfrak D}}

\allowdisplaybreaks

\title{Fractional Sobolev spaces related to an ultraparabolic operator}

\author{Antonello Pesce
\thanks{Dipartimento di Matematica, Universit\`a di Bologna, Bologna, Italy
  \textbf{E-mail}: antonello.pesce2@unibo.it} 
  \and Sascha Portaro
\thanks{Dipartimento di Matematica, Universit\`a di Bologna, Bologna, Italy
\textbf{E-mail}: sascha.portaro@unibo.it}}

\begin{document}
\maketitle

\begin{abstract}

We propose a functional framework of fractional Sobolev spaces for a class of ultra-parabolic Kolmogorov type operators satisfying the weak H\"ormander condition. 
We characterize these spaces as real interpolation of natural order intrinic Sobolev spaces recently introduced in \cite{MR4700191}, and prove continuous embeddings into $L^p$ and intrinsic H\"older spaces from \cite{MR3429628}. These embeddings naturally extend the standard  Euclidean ones, coherently with the homogeneous structure of the associated Kolmogorov group. 
Our approach to interpolation is based on approximation of intrinsically regular functions, the latter heavily relying on integral estimates of the intrinsic Taylor remainder. The embeddings exploit the aforementioned interpolation property and the corresponding embeddings of natural order intrinsic spaces.
\end{abstract}


\bigskip\noindent\textbf{Acknowledgements}: The authors are members of INDAM-GNAMPA.


\section{Introduction}\label{sec:intro}
We are interested in developing a framework of fractional Sobolev spaces for the study of evolution operators of the following type
\begin{equation}\label{PDE}
\K=\frac{1}{2}\sum_{i=1}^{d}\p_{x_i}^2-Y, \quad Y=\langle Bx,\nabla\rangle +\p_t, \qquad (t,x)\in \R\times\R^{N},
\end{equation}
where $1\le d\le N$ and $B$ is a constant $N\times N$ matrix. If $d=N$ then $\K$ is a parabolic operator while in general, for $d<N$, $\mathcal{K}$ is degenerate, and not uniformly parabolic. Nonetheless, $\K$ is a hypoelliptic operator under the H\"ormander condition
\begin{equation}\label{horm}
  \text{\rm rank Lie}(\p_{x_{1}},\dots,\p_{x_{d}},Y)=N+1.
\end{equation}
This kind of PDEs arise in many linear and non-linear models in physics: for instance, for $N=2$ and $d=1$, equation \eqref{PDE} is the {Fokker-Planck} equation for the kinetic system with Hamiltonian function 
$H(x_1,x_2)=V(x_1,x_2)+{|x_1|^2}/{2}$, where $V$ is a linear potential, and ${|x_1|^2}/{2}$ corresponds to the kinetic energy of a particle with unit mass (see \cite{Cercignani}, \cite{Desvillettes}).
In mathematical finance, equation \eqref{PDE} is related to the model used to price path-dependent contracts, such as Asian options (see, for instance, \cite{BarucciPolidoroVespri}, \cite{CIBELLI201987}, \cite{MR2791231}). For a survey of the classic theory for Kolmogorov equations see \cite{MR4069607}.

\medskip
It's well known that, under suitable assumptions on the structure of the matrix $B$, equation \eqref{PDE} endows $\R^{N+1}$ with a structure of homogeneous Lie group (see Section \ref{Lie} below), thus it is crucial to measure regularity of functions 
in terms of the group structure when dealing with differential problems related to operators of this class. In \cite{Polidoro2} a suitable notion of H\"older space was introduced, later extended in \cite{MR3429628}. Notions of $BMO$ and $VMO$ spaces have been considered in \cite{ManfrediniPolidoro}; Sobolev and Sobolev-Morrey spaces of order $k=2$ have been  introduced in \cite{BramantiCeruttiManfredini},\cite{PolidoroRagusa}, and recently extended for $k\in \N$ in \cite{MR4700191}. 

Sobolev, Hardy and Besov spaces have been studied in the context of general Lie groups since the late seventies and eighties. When the manifold is endowed with a Riemannian structure, properties of these spaces, such as embedding and interpolation theorems, as well as algebra properties are well understood: we can refer, for instance, to \cite{MR881521} and the monographs \cite{MR3024598}, \cite{MR1481970} and the references therein. 
A theory of Sobolev and Besov spaces has also been recently developed on Lie groups endowed with a sub-Riemannian structure, both for compact and non-compact manifolds: see \cite{MR3873579}, \cite{MR3944287}, \cite{MR4099628}. On the other hand, to the best of our knowledge, there are no general resources in the literature for operators that cannot be interpreted as a sub-Laplacian, such as \eqref{PDE}.

\medskip

In the remainder of this section, we introduce our notion of fractional Sobolev space $\Lambda^{n+s}_{p,q,B}$ and present the main results of the paper: Theorem \ref{th:interpolation}, which provides a characterization in terms of \textit{real interpolation} of the natural order Sobolev spaces introduced in \cite{MR4700191}, and Theorem \ref{th:embeddings}, which establishes embeddings into $L^{p}$ and the intrinsic H\"older spaces from \cite{MR3429628}. 
In the Euclidean setting these spaces correspond to a positive order Besov spaces, which reduce to Sobolev-Slobodeckij spaces when $p=q$.
In Section \ref{sec:pre} we collect some important preliminary tools and show a crucial integral estimate of the intrinsic Taylor remainder. Section \ref{sec:approx} contains the proof of our interpolation result, which is based on the approximation procedure and estimates of Proposition \ref{prop:approx}. Eventually in Section \ref{sec:emb} we prove our embeddings result.

\subsection{Intrinsic vs anisotropic H\"older and Sobolev spaces}\label{Lie}

We recall (see \cite{lanpol}) that H\"ormander's condition \eqref{horm} is equivalent to the following one: on a certain basis of $\R^N$ the matrix $B$ takes the block-form 
\begin{equation}\label{B}
  B=\begin{pmatrix}
 \ast & \ast & \cdots & \ast & \ast \\
 B_1 & \ast &\cdots& \ast & \ast \\
 0 & B_2 &\cdots& \ast& \ast \\
 \vdots & \vdots &\ddots& \vdots&\vdots \\
 0 & 0 &\cdots& B_{r}& \ast
  \end{pmatrix}
\end{equation}
where $B_j$ is a $(d_{j-1}\times d_j)$-matrix of rank $d_j$ with
\begin{equation}\label{Bd}
  d\equiv d_{0}\geq d_1\geq\dots\geq d_{r}\geq1,\qquad \sum_{j=0}^{r} d_j=N.
\end{equation}
In general, the $\ast$-blocks in \eqref{B} are arbitrary. Throughout the paper we assume the following 
\begin{assumption}\label{dilat}
$B$ is a constant $d\times d$ matrix as in \eqref{B}, where each $\ast$-block is null.
\end{assumption}
Assumption \ref{dilat} is equivalent (cf. \cite{lanpol}) to the fact that the kinetic Fokker-Planck operator \eqref{PDE} is homogeneous of degree two with respect to the family of dilations defined as follows: first of
all, consistently with the block decomposition \eqref{B} of $B$, we write $x\in\R^{N}$ as the
direct sum $x=x^{[0]}+\cdots+x^{[r]}$ where $x^{[i]}\in\R^{N}$ is defined as
  $$x_{k}^{[i]}=
  \begin{cases}
    x_{k} & \text{ if }\, \bar{d}_{i-1}<k\le\bar{d}_{i}, \\
    0 & \text{otherwise},
  \end{cases}\qquad \bar{d}_{i}:=\sum_{j=0}^{i}d_{j},\quad \bar{d}_{-1}=0,\quad i=0,\dots,r.
  $$
Then, we have $\K(u(D_{\l}))=\l^{2}(\K u)(D_{\l})$ where 
\begin{equation}\label{dilations}
  D_{\l}(t,x):=(\l^{2}t,\hD_{\l}x),\qquad \hD_{\l}x:=\sum_{i=0}^{r}\l^{2i+1}x^{[i]}.
\end{equation}
Then the $D_{\l}$-homogeneous norm on $\R^{N+1}$ is defined as
\begin{equation}\label{Bnorm}
  \|(t,x)\|_{B}=|t|^{\frac{1}{2}}+
  |x|_{B},\qquad |x|_{B}=\sum_{i=0}^{r}|x^{[i]}|^{\frac{1}{2i+1}}.
\end{equation}
\textit{Anisotropic} functional spaces on $\R^N$, defined w.r.t the component $|\cdot|_B$ of the seminorm, have been used since the seminal paper \cite{Lunardi}, when dealing with solutions of the relevant differential problem that are intended to be distributional w.r.t. the variable $t$, and {do not} benefit from the time-smoothing effect that is typical of the parabolic semigroup. {Anisotropic H\"older and Besov spaces (weighted and not) and their embeddings have been extensively exploited in many analytic applications, see for instance \cite{Lunardi}, \cite{MR4358660}, \cite{MR4554678}, \cite{MR4355925}, \cite{MR4696275}, \cite{MR4718402} and many others.}

From a functional analytic standpoint such spaces are easy to define and use, since they only require to scale the order of regularity of each set of variables $x^{[i]}$ by the appropriate factor $\frac{1}{2i+1}$, and the associated vector fields $\p_{x_1}, \dots, \p_{x_N}$ trivially commute. Therefore, the vast array of functional analytic tools available for Euclidian spaces (interpolation and duality properties, embeddings, etc.) is directly inferred for anisotropic ones.  For instance, in the case $d=1$, $N=2$, the anisotropic H\"older space of exponent $\alpha$ would consist of functions satisfying 
\begin{equation}\label{anisotropic}
|f(x)-f(\xi)|\lesssim |x_1-\xi_1|^{\alpha}+|x_2-\xi_2|^{\frac{\a}{3}}, \quad x, \xi \in \R^2.
\end{equation}
{A different approach to the construction of functional spaces on $\R^N$ 
is based on space convolutions with the associated semigroup, see for instance \cite{MR4444114}, \cite{MR4093793}.}

\medskip
If we were to include the time variable in order to fully exploit the smoothing effect of the equation in the $N+1$ variables, regularity in space and time becomes strictly intertwined. Indeed, the operator $\K$ would no longer be invariant w.r.t. Euclidean translations, but the left translations of the group $(\R^{N+1},\circ)$, where the non-commutative group law is defined by 
\begin{equation}\label{group_law}
(t,x)\circ (s,\x)=\ell_{(t,x)}(s,\x):=(s+t,\x+e^{sB}x), \qquad (t,x),\ (s,\x)\in\R^{N+1}.
\end{equation}  
Precisely, $\ell_{(t,x)}\circ\K=\K\circ\ell_{(t,x)}$ for any $(t,x)\in\R^{N+1}$.
Notice that $(\R^{N+1},\circ)$ is a group with identity $\text{Id}=(0,0)$ and inverse $(t,x)^{-1}=(-t,{-}e^{-tB}x)$, thus
the matrix $B$ uniquely identifies an \textit{homogeneous Lie group} $\cG_B:=(\R^{N+1},\circ,D(\l))$. For instance, in the previous example with 
$$B=\begin{pmatrix}
    0& 0 \\
    1 & 0
  \end{pmatrix}$$
we would have the \textit{intrinsic} H\"older space
\begin{equation}\label{intrinsic}
|f(t,x)- f(s,\x)|\lesssim \|(s,\x)^{-1}\circ (t,x)\|^\a_B\sim |t-s|^{\frac{\a}{2}}+|x_1-\x_1|^{\a}+|x_2-\x_2-(t-s)\x_1|^{\frac{\a}{3}}.
\end{equation}
Notice that for points that are on the same time level, the increments are controlled as for \eqref{anisotropic}; otherwise, as opposed to the standard parabolic case, the regularity in the Euclidean directions is somehow entangled due to the group law associated to $B$. This fact, together with the non commutativity of the fields $\p_{x_1}, \dots, \p_{x_d}, Y$ makes it difficult to build a general functional framework.

H\"older spaces defined as \eqref{intrinsic} were first used in \cite{Polidoro2}, where the author introduced the notion of solution in the \textit{Lie} sense, which regards $Y$ as a directional derivative. Intrinsic spaces associated to ${\bf d}_B$ are effective in the study of Kolmogorov equations as they comply with the asymptotic behaviour of the semigroup $\G$ of $\K$ near the pole. Indeed, it holds that 
$$\G(s,\xi;t,x)=\G(0,0; (s,\xi)^{-1}\circ (t,x)), \qquad s<t.$$
{The qualitative difference between even and odd regularity indexes makes it difficult to generalize \eqref{intrinsic} to higher orders, and for some time, only the intrinsic $C^{0,\a}$ and different notions of $C^{2,\a}$ spaces were introduced (see \cite{MR3429628}, Section 3.2 for a review of the intrinsic H\"older spaces in the literature). Eventually in \cite{MR3429628} a general definition of $C^{n,\a}$ was established: owning to an estimate of the intrinsic Taylor remainder, it turns out that it is equivalent to set a notion of regularity by means of the intrinsic distance
or to set the appropriate regularities with respect to the H\"ormander vector fields $\{\p_{x_1},\dots,\p_{x_d},Y\}$, according to their \textit{formal degree}. For instance, condition \eqref{intrinsic} is equivalent to set $u\in C^{\a}$ along the direction $x_1$ and $u\in C^{\a/2}$ along the integral curves of the vector field $Y$. In the following paragraphs, we recall the general intrinsic H\"older spaces which are relevant for our embeddings, and introduce the fractional Sobolev spaces.}

\medskip
Let $h\mapsto e^{h X }z$ denote the integral curve of a Lipschitz vector field $X$ starting from
$z\in\R^N$, defined as the unique solution of
\begin{equation}
\begin{cases}
 \frac{d}{dh}e^{h X }z= X\left(e^{h X }z\right),\qquad  &h\in\R, \\
 e^{h X }z\vert_{h=0}= z.
\end{cases}
\end{equation}
For the vector fields in $\K$, we have
\begin{equation}\label{eq:def_curva_integrale_campo}
 e^{h \p_{x_{i}} }(t,x)=(t,x+h \mathbf{e}_i),\qquad 
 e^{h Y }(t,x)=(t+h,e^{h B}x),
\end{equation}
where $\mathbf{e}_i$ is the $i$-th element of the canonical basis of $\R^{N}$. A simple computation shows that, for any $z,\z\in \R^{N+1}$,
\begin{equation}\label{invariance_law}
 \z^{-1}\circ e^{\delta Y}z=e^{\delta Y}(\z^{-1}\circ z), \qquad \z^{-1}\circ e^{\delta \p_{x_i}}z=e^{\delta \p_{x_i}}(\z^{-1}\circ z),
 \qquad i=1,\dots , d.
\end{equation}
Moreover, we have  (see, for instance, \cite{MR4694350})
\begin{equation}\label{homogeneity}
 D_\l e^{\delta Y}(z)=e^{\delta\lambda^2Y}\left(D_\l z\right), \qquad D_\l e^{\delta
 \p_{x_i}}(z)=e^{\delta \l^{2j+1} \p_{x_i}}D_\l z, \qquad i=\bar{d}_{j-1}+1,\dots ,\bar{d}_{j}.
\end{equation}
\begin{definition}
Let $m_{X}$ be a formal weight associated to the vector field $X$. For $\a\in\,]0,m_{X}]$, we say
that $u\in C_{X}^{\alpha}$ if the quasi-norm
\begin{equation}\label{CaX}
 \|u\|_{C^{\a}_{X}}:=
 \sup_{z\in\R^N\atop h\in\R\setminus\{0\}} \frac{\left|u\left(e^{h X }z\right)-u(z)\right|}{|h|^{\frac{\alpha}{m_{X}}}}
\end{equation}
is finite.
\end{definition}

Hereafter, we set {\it the formal weight of the vector fields $\p_{x_{1}},\dots,\p_{x_{d}}$ equal
to one and the formal weight of $Y$ equal to two}, which is coherent with the homogeneity of the
Fokker-Planck operator $\K$ with respect to the dilations $D_{\l}$ in \eqref{dilations}. From
\cite{MR3429628} we recall the following
\begin{definition}[\bf Intrinsic H\"older spaces]\label{def:C_alpha_spaces}
For $\a\in\,]0,1]$ we define the H\"older quasi-norms 
\begin{align}
  \norm{u}_{C^{0,\a}_{B}}&:=\sup_{\R^{N+1}}|u|+\sum_{i=1}^{d}\norm{u}_{C^{\a}_{\partial_{x_i}}}+\norm{u}_{C^{\a}_{Y}},\\
  \norm{u}_{C^{1,\a}_{B}}&:=\sup_{\R^{N+1}}|u|+\norm{\nabla_{d}u}_{C^{0,\a}_{B}}+\norm{u}_{C^{\a+1}_{Y}},
\intertext{where $\nabla_{d}:=(\p_{x_{1}},\dots,\p_{x_{d}})$ and inductively, for $n\ge 2$,}
  \norm{u}_{C^{n,\a}_{B}}&:=\sup_{\R^{N+1}}|u|+\norm{\nabla_{d}u}_{C^{n-1,\a}_{B}}+\norm{Y u}_{C^{n-2,\a}_{B}}.
\end{align}
\end{definition}
\medskip
Following this approach, it is natural to introduce a notion of \textit{intrinsic Besov spaces} which are essentially built by replacing the sup norms in \eqref{CaX} with $L^p$ and $L^q$ norms w.r.t. $z$ and $h$ respectively, and is {coherent} with the work \cite{MR1762582} on fractional Sobolev-Slobodeckij spaces of order $s\in (0,1)$ associated to a H\"ormander structure. More precisely, for $p,q\ge 1$, we set the fractional Besov quasi-norm of order $s$, $0<s<m_X$ along a vector field $X$
\begin{equation}
{[u]_{X,s,p,q}:=\begin{cases}
\left( \int_{ \R}\|u(e^{hX}\cdot)-u\|^q_{{p}}\frac{dh}{|h|^{1+q\frac{s}{m_X}}}\right)^{\frac 1q}, & q\in [1,\infty)\\
\sup_{h\neq 0}\frac{\|u(e^{hX}\cdot)-u\|_{{p}}}{|h|^{\frac{s}{m_X}}}, & q=\infty.
\end{cases}}
\end{equation}
Then for $p,q\ge 1$ and $s\in (0,1)$ we set
\begin{align}
  |u|_{s,p,q}&:=\sum_{i=1}^d [u]_{\p_{x_i},s,p,q}+[u]_{Y,s,p,q},\\ \label{W2pp2}
  |u|_{1+s,p,q}&:=|\nabla_{d}u|_{s,p,q}+[u]_{Y,1+s,p,q},
\intertext{and inductively, for $n\ge 2$,}\label{Wnppn}
  |u|_{n+s,p,q}&:=|\nabla_{d}u|_{n+s-1,p,q}+|Yu|_{n+s-2,p,q}, \quad n\in \N, \; n\ge 2
\end{align}
{Hereafter we will denote $\N_0:=\N\cup \{0\}$}.
\begin{definition}[\bf Intrinsic Besov spaces]\label{def1}
For $p,q\ge 1$, $n\in \N_0$ and $s\in (0,1)$, {we define the Besov $\La^{n+s}_{p,q,B}$ as the space of measurable functions on $\R^{N+1}$ such that the following quasi-norm is finite:}
\begin{align}\label{eq:Wnpq}
\|u\|_{\La^{n+s}_{p,q,B}}&:=\|u\|_p+|u|_{n+s,p,q}.
\end{align}
\end{definition}

\subsection{Main results}
\subsubsection{Interpolation characterization}
In the Euclidean setting, it is well known that fractional Sobolev-Slobodeckij spaces and positive order Besov spaces can be characterized by \textit{real interpolation} of natural order Sobolev spaces, see for instance the references \cite{MR3024598}, \cite{MR2424078}. We recover the analogous result where the natural order Euclidean Sobolev spaces are replaced by the intrinsic Sobolev spaces introduced in \cite{MR4700191}, which are defined as follows: for any  
$p\ge 1$ set
\begin{align}
  |u|_{1,p}&:=\sum_{i=1}^d \|\p_{x_i}u\|_{p}+[u]_{Y,\frac{1}{2},p,p},\\ 
  |u|_{2,p}&:=|\nabla_{d}u|_{1,p}+\|Yu\|_p,
\intertext{and inductively, for $n\ge 2$,}
  |u|_{n,p,q}&:=|\nabla_{d}u|_{n-1,p}+|Yu|_{n-2,p}, \quad n\in \N, \; n\ge 3.
\end{align}
\begin{definition}[\bf Intrinsic Sobolev spaces]
For $p\ge 1$ and $n\in \N_0$ and $s\in (0,1)$, {the Sobolev quasi-norm is defined as}
\begin{align}\label{eq:Wnp}
\|u\|_{W^{n,p}_{B}}&:=\|u\|_p+|u|_{n,p}.
\end{align}
\end{definition}

Before we state our interpolation result, let us first recall what a real interpolation space is: for a comprehensive
presentation of the subject we refer, for instance, to \cite{MR3726909}, \cite{MR2328004} and
\cite{MR2424078}. Given two real Banach spaces $Z_1, Z_2$, we write $Z_1=Z_2$ if $Z_1$ and $Z_2$ have the same
elements with equivalent norms and we write $Z_1\subseteq Z_2$ if $Z_1$ is continuously embedded in
$Z_2$. The pair $(Z_1, Z_2)$ is called an \textit{interpolation pair} if both $Z_1$ and $Z_2$ are
continuously embedded in some Hausdorff topological vector
space: in this case, the intersection $Z_1\cap Z_2$ and the sum 
$Z_1+Z_2$ endowed with the norms
  $$\|u\|_{Z_1\cap Z_2}:=\max\{\|u\|_{Z_1},\|u\|_{Z_2}\}, \qquad \|u\|_{Z_1+ Z_2}:=
  \inf_{u_1\in Z_1,\, u_2\in Z_2\atop u=u_{1}+u_{2}}(\|u_1\|_{Z_1}+\|u_2\|_{Z_2}),$$
are Banach spaces. For any $\l>0$ and $u\in Z_{1}+Z_{2}$, we set
\begin{equation}\label{ae26}
 K(\l,u)\equiv K(\l,u;Z_1,Z_2):=
  \inf_{u_1\in Z_1,\, u_2\in Z_2\atop u=u_{1}+u_{2}}(\|u_1\|_{Z_1}+\l\|u_2\|_{Z_2}).
\end{equation}
Any Banach space $E$ such that  
\begin{equation}\label{intermediate}
Z_1\cap Z_2\subseteq E\subseteq Z_1+Z_2,
\end{equation}
is called an \textit{intermediate space}. Among these, for $0<\theta<1$ and $1\le p\le \infty$, we
have the \textit{real interpolation space} $(Z_1,Z_2)_{\theta,p}$ consisting of $u\in Z_{1}+Z_{2}$
such that
\begin{equation}\label{k_interpolation}
  \|u\|_{\theta,p}:=\|\l\mapsto \l^{-\theta}K(\l,u)\|_{L^{p}(\R_{>0},\frac{d\l}{\l})}<\infty.
\end{equation}

We are ready to state our first result.
\begin{theorem}\label{th:interpolation}
For any $n\in \N_0$, $s\in (0,1)$ and {$p,q\ge 1$}, we have
\begin{align*}
\Lambda^{n+s}_{p,q,B}&=\big(L^p,W^{n+1,p}_{B}\big)_{\frac{n+s}{n+1},q}.
\end{align*}
\end{theorem}

\medskip
\noindent\emph{Intermediate spaces and reiteration}: the set of intermediate spaces $E$ for the pair $(Z_1,Z_2)$ such that 
$$(Z_1,Z_2)_{\theta,1}\subseteq E\subseteq (Z_1,Z_2)_{\theta,\infty}$$
for some $\theta\in (0,1)$ constitutes a special class of intermediate spaces denoted by $\mathcal{H}_{\theta}(Z_1,Z_2)$; for instance $(Z_1,Z_2)_{\theta,p}\in \mathcal{H}_{\theta}(Z_1,Z_2)$. If $E_i\in \mathcal{H}_{\theta_i}(Z_1,Z_2)$, $i=1,2$ {then} it holds that 
$$(Z_1,Z_2)_{(1-\l)\theta_1+\l\theta_2,p}=(E_1,E_2)_{\l,p}, \quad p\ge 1$$
with equivalence of norms (see Theorem 1.23 in \cite{MR3753604}). This result is known as the \textit{Reiteration Theorem} and it is useful to obtain equivalent characterizations for interpolation spaces. 

Now by Theorem 5.3 in \cite{MR4700191} we have that $W^{n,p}_B\in  \mathcal{H}_{\frac{n}{m}}(L^{p},W^{m,p}_B)$ for any $n,m\in \N$, $m>n$. Together with the Reiteration Theorem, this allows to infer the following equivalences: we have 
%
%
\begin{align*}
\Lambda^{n+s}_{p,q,B}&=\big(L^p,W^{n+1,p}_{B}\big)_{\frac{n+s}{n+1},q}
\intertext{(by reiteration, for $0<s\leq s_1$, $n_1\in \N_0$, $n_1\geq n$ and $q_1\in [1,\infty]$)}
&=\big(L^p,\Lambda^{n_1+s_1}_{p,q_1,B}\big)_{\frac{n+s}{n_1+s_1},q}
 \intertext{(by reiteration, for $n_1, n_2\in \N_0$, $s_1, s_2 \in (0,1)$ s.t. $n_1+s_1\le n+s\le n_2+s_2$ and $q_1,q_2\in [1,\infty]$)}
 &=\big(\La^{n_1+s_1}_{p,q_1,B},\La^{n_2+s_2}_{p,q_1,B}\big)_{\theta,q}, \qquad n+s=(1-\theta)(n_1+s_1)+\theta (n_2+s_2).
\end{align*}

\subsubsection{Embeddings into Lorentz and intrinsic H\"older spaces}

With reference to \eqref{Bd}, we set
\begin{equation}\label{homdim}
  \Q:= 2+\sum_{k=0}^{r}(2k+1)d_{k}
\end{equation}
the \textit{homogeneous dimension} of $\R^{N+1}$ with respect to $D_{\l}$. Then we have the following
%

\begin{theorem}\label{th:embeddings}Let $n\in \N_0$, $s\in (0,1)$, {$p\in (1,\infty)$ and $q\ge 1$}. The following inclusions hold with continuity:
\begin{itemize}
\item[i)]If $(n+s)p={\bf d}$, then $\La^{n+s}_{p,1,B}\subseteq L^{\infty}$;
\item[ii)]If $(n+s)p<{\bf d}$, 
then $\La^{n+s}_{p,q,B}\subseteq L^{p',q}$, with $\frac{1}{p'}=\frac{1}{p}-\frac{n+s}{\bf d}$;
\item[iii)] If $(n+s)p>{\bf d}$, then $\La^{n+s}_{p,q,B}\subseteq C^{n-k,k+s-\frac{\bf d}{p}}_{B}$, 
\end{itemize}
where $k$ is the smallest integer s.t. $(k+s)p>{\bf d}$. 
\end{theorem}

{Theorem} \ref{th:embeddings} gives the expected embeddings, coherently with the homogeneous structure of the group. Indeed these are analogous to the usual Euclidean ones, once we replace $N+1$ with {the homogeneous dimension of $N+1$, and the Euclidean H\"older spaces} with the intrinsic ones. Theorems \ref{th:interpolation} and \ref{th:embeddings} show that our functional framework is compatible with the notion of interpolation and also with the H\"older functional framework introduced in \cite{MR3429628}.

%

\section{Intrinsic weak derivatives and Taylor expansion}\label{sec:pre}
\subsection{Preliminaries}
Here we gather some properties that {are} useful in the next sections. 

\medskip
Hereafter, in the context of our proofs {we often use} the notation $A \lesssim B$, meaning that $A \le cB$ for some
positive constant $c$ which may depend on the quantities specified in the corresponding statement. {We also exclusively write} {the proofs for $q<\infty$, the case $q=\infty$ being an easier modification}.

\begin{lemma}\label{lem0}[{\cite{MR1751429}, Proposition $5.1$}] There exists $m=m(B)\ge 1$ such that
\begin{equation}\label{normcontrol}
 \|\z^{-1}\circ z\|_B\le m(\|\z\|_B+\|z\|_B)\quad m^{-1}\|z\|_B\le\|z^{-1}\|_B\le  m\|z\|_B \quad
 z,\z\in \R^{N+1}.
\end{equation}
\end{lemma}
Note that, since $e^{\delta Y}z=z\circ (\delta,0)$, by \eqref{normcontrol} we have
\begin{equation}\label{normcontrol2}
  \frac{1-mc}{c}\|z\|_{B}\le \|e^{\delta Y}z\|_B\le m(1+c)\|z\|_B,
\end{equation}
for any $|\d|^{\frac 12}\le c\|z\|_B$ with $c\in\,]0,\frac{1}{m}[$.

\begin{remark}(see, for instance \cite{MR4694350})\label{r1} 
The matrix $B$ is nilpotent of degree $r+1$. Moreover, the exponential matrix
$$e^{\delta B}=I_N+\sum_{j=1}^r\frac{B^j}{j!}\delta^j$$
is a lower triangular with diagonal $(1,\dots, 1)$ and therefore it has determinant equal
to 1.
\end{remark}

The next result shows that the Besov quasi-norms scale homogeneously w.r.t. derivatives of the same intrinsic order. 
\begin{lemma}
For any $n\in \N_0$, $s\in (0,1)$, $p,q\ge 1$ and $\l>0$ we have
\begin{equation}\label{scaling}
 \|u( D_\l)\|_p=\l^{-\frac{\Q}{p}}\|u\|_p,\qquad |u( D_\l)|_{n+s,p,q}=\l^{n+s-\frac{\Q}{p}}|u|_{n+s,p,q}
\end{equation}
\end{lemma}
\proof The first equality follows by a simple change of variable. Next, let $X\in \{\p_{x_1},\dots, \p_{x_d},Y\}$. We have 
\begin{align*}
[u(D_{\l})]^q_{X,s,p,q}&=\int_{\R}\|u(D_{\l}e^{hX}{\cdot})-u(D_{\l})\|^q_{{p}}\frac{dh}{|h|^{1+\frac{s}{m_X}q}}\\
&=\l^{m_X+sq}\int_{\R}\left(\int_{\R^{N+1}}|u(e^{h\l^{m_X}}D_{\l})-u(D_{\l})|^pdz\right)^{\frac{q}{p}}\frac{dh}{|\l^{m_X} h|^{1+\frac{s}{m_X}q}}
\intertext{(by the change of variable $z'=D_{\l}z$)}
&=\l^{m_X+{q(s-\frac{\Q}{p})}}
\int_{\R}\|u(e^{h\l^{m_X} X}\cdot)-u\|^q_{{p}}\frac{dh}{|\l^{m_X} h|^{1+\frac{s}{m_X}q}}
\intertext{(by the change of variable $h'=\l^{m_X}h$)}
&=\l^{q(s-\frac{\Q}{p})}[u]_{X,s,p,q}^q.
\end{align*}
This gives \eqref{scaling} for $n=0$. For higher orders it suffices to notice that $Xu(D_{\l}z)=\l^{m_X}(Xu)(D_{\l}z)$ by \eqref{homogeneity} and proceed by induction.
\endproof 

For any $\b\in \N_0^{N}$ and denote
\begin{equation}\label{B_lenght}
 \p^\b=\p^{\b_1}_{x_1}\cdots \p^{\b_N}_{x_N}, \qquad 
 {\bbB=\sum_{i=0}^r(2i+1)\sum_{k=1+\bar{d}_{i-1}}^{\bar{d}_{i}}\b_{k}.}
\end{equation}
The following result is analogous to Lemma 3.1 and Proposition 3.2 in \cite{MR4700191}. It shows that even if the Besov quasi-norms only control weak derivatives made of compositions of $\p_{x_1},\dots,\p_{x_d}$ and $Y$, a distribution $u\in \La^{n+s}_{p,q,B}$ can actually {support} all the weak derivatives
of intrinsic order $l$, namely $Y^k\p^\b u$ with $2k+\bbB=l$, if $n$ is high enough.

\begin{proposition}\label{prop:prop1}
Let $u\in \La^{n+s}_{p,q,B}$, $\b\in \N_0^{N}$ and $k\in \N_0$. If $n\ge 2k+\bbB$, then 
\begin{equation}\label{intrinsic_D}
 Y^k\p^\b u\in \La^{n+s-2k-\bbB}_{p,q,B}.
\end{equation}
\end{proposition}

\subsection{Taylor inequality}\label{sec:Taylor}

The following crucial estimate holds.
\begin{theorem}\label{th:Taylor}
Let $n\in \N_0$, $s\in (0,1)$ and $p,q\ge 1$. There exists a constant $c=c(n,s,p,q,B)$ such that, for any $u\in  \La^{n+s}_{p,q,B}$ we have 
\begin{equation}\label{eq:Taylor0}
\int_{\R^{N+1}}\|u-T_nu(\cdot\circ \z,\cdot)\|^q_{{p}}\frac{d\z}{\|\z\|_{\bf B}^{{\bf d}+{(n+s)}q}} 
\le c {|u|_{n+s,p,q}},
\end{equation}
where  $T_nu(\z,\cdot)$ is the $n$-th order $B$-Taylor polynomial (see \cite{MR3429628}) of $u$ around
$\z=(\t,\x)$, formally defined as
\begin{equation}\label{eq:def_Tayolor_n}
 T_n u(\z,z):= \sum_{{0\leq 2 k + |\beta|_B \leq n}}\frac{(t-\t)^k (x-e^{(t-\t)B}\xi)^{\beta}}{k!\,\beta!}
 Y^k \partial_{\xi}^{\beta}u(\t,\xi),\qquad
 z=(t,x)\in\R^{N+1}.
\end{equation}
\end{theorem}

\begin{remark}Notice that, for $n=0$ and $p=q$, the LHS is the usual Sobolev-Slobodeckij seminorm for $W^{s,p}$, where the Euclidean translation and distance are replaced by the intrinsic translation and homogeneous norm of the group.
\end{remark}

\proof[Strategy of the proof] The general proof of Theorem \ref{th:Taylor} is a non-trivial adaptation of the two-scale induction procedure devoloped in \cite{MR3429628} for the (intrinsic) H\"older version ($p=q=\infty$) and the proof of 
{Theorem 4.1 in \cite{MR4700191}}. The point is to provide controls for the increments between the integration points $z$ in the inner norm, and $z\circ \z$, by moving along the integral curves of the H\"ormander vector fields $\p_{x_1},\dots,\p_{x_d}, Y$ and their commutators, the latter being approximated using classical techniques from control theory. 
To show how to modify the proof of Theorem in 4.1 in \cite{MR4700191} and highlight the main novelty of the fractional case, we present the proof in the simpler case of the Langevin operator in $\R^3$, $\mathcal{K}=\frac{1}{2}\p^2_{x_1}-x_1\p_{x_2}-\p_t$ (${\bf d}=6$). The generalization to the full Kolmogorov chain and the induction procedure is analogous to \cite{MR4700191}.

\medskip
Let $\z=(\t,\x)\in \R\times \R^2$ and $z=(t,x)\in\R\times \R^2$. From \eqref{group_law}, we have in particular
$$z\circ \z=(t+\t,\x+e^{\t B}x)=(t+\t,x_1+\x_1,x_2+\x_2+\t x_1).$$ 
To reach the point $z\circ \z$ starting from $z$ moving along the H\"ormander vector fields, we first connect $z$ to 
$z_{-1}=e^{\t Y}z=(t+\t,x_1,x_2+\t x_1)$, then connect $z_{-1}$ to $z_0=e^{\x_1\p_{x_1}}z_{-1}=(t+\t ,x_1+\x_1,x_2+\t x_1)$. At this point it only remains to move by $\x_2$ in direction $\p_{x_2}$, that is the commutator $[\p_{x_1},Y]$. 
This is done by following the integral paths of the vector fields $\p_{x_1}$  and $Y$ as follows:
$$z\circ \z=e^{-\x_2^{\frac{2}{3}}Y}\left(e^{-\x_2^{\frac{1}{3}}\p_{x_1}}\left(e^{\x_2^{\frac{2}{3}}Y}\left( e^{\x_2^{\frac{1}{3}}\p_{x_1}}z_0\right)\right)\right).$$
Next we separately control the increments, thanks to Minkowski's inequality and the fact
that the Jacobian of the changes of variables $z\mapsto e^{\delta \p_{x_1}}z$ and $z\mapsto e^{\delta Y}z$ have determinant equal to one (see Remark \ref{r1}). 
For a permutation $\{i,j\}$ of $\{1,2\}$ we write
\begin{align}
&\int_{\R^{3}} \left\|u(e^{\x_i^{\frac{1}{2i-1}}\p_{x_1}}\cdot)-u\right\|^q_{{p}} 
\frac{d\z}{\|\z\|_{\bf B}^{6 + s q}} \nonumber \\ 
    &\quad=\int_{\R^{3}} {\left\|u(e^{\x_i^{\frac{1}{2i-1}}\p_{x_1}}\cdot)-u\right\|^q_{{p}}}
    \frac{d\t  d\x_i d\x_j}{|\x_i|^{\frac{6 + s q}{2i-1}} \left (1+\frac{|\t |^{\frac{1}{2}}+|\x_j|^{\frac{1}{2j-1}}}{|\x_i|^{\frac{1}{2i-1}}} \right )^{6+s q}}
\intertext{(by the change of variables $\x_i=h^{2i-1}$, $\x_j=h^{2j-1}\bar{\x}_j$, $\t =h^2\bar{\t }$)}
&\quad=\int_{\R}\|u(e^{h\p_{x_1}}\cdot)-u\|^q_{p}\frac{dh}{|h|^{1+sq}}
\int_{\R^2}\frac{{(2i-1)} d\bar{\t } d\bar{\x}_j}{(1+|\bar \t |^{\frac{1}{2}}+|\bar \x_j|^{\frac{1}{2j-1}})^{6+sq}}
{= c_{q} [u]_{\p_{x_1},s,p,q}^q.} \label{eq:T1}
\end{align}
{Similarly, by the change of variables 
$\xi_1 = |\t |^{\frac{1}{2}} \bar{\xi}_1$, $\xi_2 = |\t |^{\frac{3}{2}} {\bar{\xi_2 }}$, we also have
\begin{align}
&\int_{\R^{3}} \left\|u(e^{sY}\cdot)-u\right\|^q_{{p}} \frac{d\z}{\|\z\|_{\bf B}^{6 + sq}}
=([u]_{Y,s,p,q})^q\int_{\R^2}\frac{d\bar{\x}_1d\bar{\x}_2}{(1+|\bar{\x}_1|+|\bar{\x}_2|^{\frac 13})^{6+sq}}=c'_q[u]_{Y,s,p,q}^q.\label{eq:T2}
\end{align}
{Lastly, by the change of variables 
${\x_2=\sgn(h)|h|^{\frac{3}{2}}}$, ${\x_1=|h|^{\frac{1}{2}}\bar{\x}_1}$, $\t =|h|\bar{\t }$, we have}
\begin{align}
&\int_{\R^{3}} \left\|u(e^{{\x_2^{\frac{2}{3}}}Y}\cdot)-u\right\|^q_{{p}}\frac{d\z}{\|\z\|_{\bf B}^{6 + sq}}
\nonumber\\
    &\quad {=\int_{\R}\|u(e^{|h|Y}\cdot)-u\|^q_{p}\frac{dh}{|h|^{1+{\frac{s}{2}}q}}}
\int_{\R^2}\frac{{3} d\bar{\t } d\bar{\x}_{{1}}}{2(1+|\bar \t |^{\frac{1}{2}}+|{\bar \x_1}|)^{6+sq}}
 {\le} \frac{c_q}{2}[u]_{Y,s,p,q}^q.\label{eq:T3}
\end{align}
Gathering together \eqref{eq:T1}, \eqref{eq:T2}, \eqref{eq:T3}, we eventually infer estimate \eqref{eq:Taylor0} for $n=0$.
\endproof

\medskip
\noindent Importantly, Theorem \ref{th:Taylor} provides the exact fractional regularity in any Euclidean direction. 
For instance, let $u\in \Lambda^s_{p,q,B}$ and $i\in \{\bar d_{k-1}+1,\dots, \bar d_k\}$: let us consider the fractional seminorm $[u]_{\p_{x_i},\g,p,q}$ for some $\g\in (0,1)$. Coherently with \eqref{dilations}, this should be finite for {$\g=\frac{s}{2k+1}$}.

\medskip
For any point $y\in \R^{N}$, denote 
$\hat y=(y_1,\dots, y_{i-1},y_{i+1},\dots,y_N)\in \R^{N-1}$ and 
$$\|(\rho,\hat y)\|_{B,i}=|\rho|^{\frac 12}+\sum_{j=1\atop j\neq i}^n|y_j|^{\frac{1}{\sigma_j}},$$
 where $\sigma_j\geq 1$, $j=1,\dots N$ are such that 
$\hat D_{\l}y=(\l^{\sigma_1}y_1,\dots,\l^{\sigma_N}y_N)$. By the change of variables 
$$ h^{2k+1}=\x_i, \quad  h^{2}{\bar s}=s, \quad  h^{2l+1}\bar {\x_j}=\x_j, \quad  \text{for any} \quad j\in  \{\bar d_{l-1}+1,\dots, \bar d_l \}$$
we have 
\begin{align*}
&\int_{\R^{N+1}}\|u(e^{\x_i \p_{x_i}}\cdot)-u\|_{p}^q\frac{d\z}{\|\z\|_B^{{\bf d}+{sq}}}\\
&\quad = (2k+1)\int_{\R}\|u(e^{h^{2k+1} \p_{x_i}}\cdot)-u\|_{p}^q\frac{dh}{|h|^{1+sq}} \underbrace{\int_{\R^N}\frac{d\bar s d\hat{\bar\x}}{(1+\|(\bar s,\hat{\bar\x})\|_{B,i})^{{\bf d}+sq}}}_{=:c_{i,q}}
\intertext{(by the change of variables $\bar h=h^{2k+1}$)}
&\quad =  c_{i,q}\int_{\R}\|u(e^{\bar h \p_{x_i}}\cdot)-u\|_{p}^q\frac{dh}{|h|^{1+\frac{s}{2k+1}q}}=c_{i,q}[u]_{\p_{x_i},\frac{s}{2k+1},p,q}^q.
\end{align*}
Therefore we get the expected control
\begin{align*}
[u]_{\p_{x_i},\frac{s}{2k+1},p,q}^q&=c_{i,q}^{-1}\int_{\R^{N+1}}\|u(e^{\x_i \p_{x_i}}\cdot)-u\|_{p}^q\frac{d\z}{\|\z\|_B^{{\bf d}+{sq}}}\\
&\lesssim \int_{\R^{N+1}}\left(\int_{\R^{N+1}}|u(z)-u(z\circ \z)|^pdz\right)^{\frac{q}{p}}\frac{d\z}{\|\z\|_B^{{\bf d}+{sq}}} \lesssim
{|u|_{s,p,q}^q}<\infty.
\end{align*}
 
More generally, we have the following

\begin{corollary}
Let $n\in \N_0$, $s\in (0,1)$, $k\in\N_{0}$, $\b\in \N_0^N$ and $j\in \N$ such that $2j+1>n+s-2k-\bbB\ge 0$.
Then there exists $c=c(n,s,k,\b,j)$ such that, 
for every $u\in \La^{n+s}_{p,q,B}$ and $i\in \{\bar d_{j-1}+1,\dots, \bar d_j\}$ we have
 $$[Y^k\p^{\b}u]_{\p_{x_i},\frac{n+s-2k-\bbB}{2j+1},p,q} \le c {|u|_{n+s,p,q}}.$$
\end{corollary}

\section{Approximation and proof of the main interpolation result}\label{sec:approx}
Following \cite{MR4694350}, we define the smooth approximation of $u\in \La^{n+s}_{p,q,B}$ as 
\begin{align}\label{approx}
  u_{n,\eps}(z):=
  \int_{\R^{N+1}}T_{n}u(\z,z)\phi\left(D_{\eps^{-1}}(\z^{-1}\circ
  z)\right)\frac{d\z}{\eps^{\Q}},\qquad \e>0,
\end{align}
where $\phi$ is a test function supported on $\|z\|_B\le 1$. In particular
\begin{equation}\label{eq:mollifier}
 \int_{\R^{N+1}}\phi\left(D_{\eps^{-1}}(\z^{-1}\circ z)\right)\frac{d\z}{\eps^{\Q}}=\int_{\|\z\|_B\le 1}\phi(\z)d\z=1.
\end{equation} Then we have the following 

\begin{proposition}[Approximation]\label{prop:approx}
Let $s\in (0,1)$, $n\in \N_0$ and $m\in\N$, $m>n$. Then there exists a constant 
$c=c(s,n,m,{\bf d})$ such that, for any $u\in \La^{n+s}_{p,q}$ we have
\begin{align}
\left\|\eps\mapsto \eps^{-(n+s)}\|u-u_{n,\eps} \|_{L^p} \right\|_{L^q({\R_{>0}},\frac{d\eps}{\eps})}
&\le c {|u|_{n+s,p,q}}\label{eq:approx1} \\
\left\|\eps\mapsto \eps^{m-(n+s)}\|u_{n,\eps} \|_{W^{m,p}_B} \right\|_{L^q({\R_{>0}},\frac{d\eps}{\eps})}
&\le c {|u|_{n+s,p,q}}.\label{eq:approx2}
\end{align}
\end{proposition}

A key tool in the proof is the following property(see Lemma 5.1 in \cite{MR4700191}).
\begin{lemma}\label{LEM}
For any $u\in W^{n,p}_B$ and ${z,\z}\in \R^{N+1}$ we have
\begin{align}
 \p_{x_i}T_nu(\z,z)&=T_{n-1}(\p_{i}u)(\z,z), \qquad n\ge 1, \ i=1,\dots,d, \label{eq:Taylor1}\\
 Y_zT_nu(\z,z)&=T_{n-2}(Yu)(\z,z), \qquad n\ge 2. \label{eq:Taylor2}
\end{align}
\end{lemma}

\proof[Proof of Proposition \ref{prop:approx}] We organize the proof in 3 main parts. 

\medskip 
\noindent \textit{PART 1: Proof of \eqref{eq:approx1}}.
By  \eqref{approx} and the change of variable $\eta=\z^{-1}\circ z$ we have
\begin{align*}
&\left\|\eps\mapsto \eps^{-(n+s)}\|u-u_{n,\eps} \|_{L^p} \right\|^q_{L^q({\R_{>0}},\frac{d\eps}{\eps})}\\
    &=\int_0^{{\infty}}{\eps^{-(n+s)q} \left ( \int_{\R^{N+1}} \left | \int_{\R^{N+1}}(T_{n} u(\z,z) - u(z)) 
    \phi \left(D_{\eps^{-1}}(\z^{-1}\circ z)\right)\frac{d\z}{\eps^{\Q}} \right |^p dz \right )^{\frac{q}{p}} \frac{d\eps}{\eps}} \\
    &\lesssim \int_0^{{\infty}}{\eps^{-(n+s)q} \left ( \int_{\R^{N+1}} \left | \int_{\R^{N+1}}(T_{n} u(z \circ \eta^{-1},z) - u(z)) 
    \phi \left(D_{\eps^{-1}} \eta \right)\frac{d\eta}{\eps^{\Q}} \right |^p dz \right )^{\frac{q}{p}} \frac{d\eps}{\eps}} \\
    \intertext{(by Minkowski inequality)}
    &= \int_0^{{\infty}}{\eps^{-(n+s)q} \left ( \int_{\R^{N+1}} \|T_{n} u(\cdot \circ \eta^{-1},\cdot) - u\|_p \phi \left(D_{\eps^{-1}} \eta \right) \frac{d\eta}{\eps^{\Q}} \right )^q \frac{d\eps}{\eps}} \\  
  \intertext{(by H\"older inequality and \eqref{eq:mollifier})}  
    &\leq \int_{\| \eta \|_{\bf B}}^{{\infty}}{\eps^{-(n+s)q} \int_{\R^{N+1}} \|T_{n} u(\cdot \circ \eta^{-1},\cdot) - u\|_p^q \phi \left(D_{\eps^{-1}} \eta \right) \frac{d\eta}{\eps^{\Q}} \frac{d\eps}{\eps}} \\
    &\lesssim \int_{\R^{N+1}} \|T_{n} u(\cdot \circ \eta^{-1},\cdot) - u\|_p^q \frac{d\eta}{\| \eta \|_{B}^{\Q+(n+s)q}} 
    \lesssim {|u|_{n+s,p,q}^q}
\end{align*}
by Theorem \ref{th:Taylor}. This concludes the proof of \eqref{eq:approx1}.

\medskip 
\noindent \textit{PART 2: Partial proof of \eqref{eq:approx2}}.
Next we check \eqref{eq:approx2} provided the preliminary estimates \eqref{eq:estim1_b}-\eqref{eq:estim2_b} below hold.
Recalling Corollary 3.4 in \cite{MR4700191}, the quasi-norm $|u_{n,\eps}|_{m,p,B}$ is equivalent to
$$\sum_{2k+\bbB=m}\|Y^k\p^\b u_{n,\eps}\|_p+\sum_{2k+\bbB=m-1}[Y^k\p^\b u_{n,\eps}]_{Y,\frac 12,p}.$$
We denote by $\DB^l$ any weak derivative of intrinsic order $l$, that is
$\DB^l=Y^k\p_{x}^\b$ for some $k\in \N_0$, $\b \in \N_0^N$ such that $2k+\bbB=l$, and denote
\begin{align}
 I^{(n,l)}_\eps u(z):=\int_{\R^{N+1}}(T_{n}u(\z,z)-u(z)){\DB^l\phi}\left(D_{\eps^{-1}}(\z^{-1}\circ z)\right)\frac{d\z}{\eps^{\Q}}.
\end{align}
Assume for a moment that, for any $i\in 0,\dots,n$ we have 
\begin{align}
    \left \| \eps \mapsto {\eps^{m-(n+s)}}\|I^{(n-i,{m-i})}_{\eps} \DB^i u\|_p \right\|_{L^q({\R_{>0}},d\eps/\eps)} 
    &\lesssim {| \DB^i u|_{n-i+s,p,q}}; \label{eq:estim1_b} \\ 
\left \| \eps \mapsto {\eps^{m-(n+s)}}[I^{(n-i,{m-i-1})}_{\eps} \DB^i u]_{Y,1,p} \right\|_{L^q({\R_{>0}},d\eps/{\eps})} 
&\lesssim {|\DB^i u|_{n-i+s,p,q}}. \label{eq:estim2_b}
\end{align}
Then, by Lemma \ref{LEM} and the Leibniz rule, we have
\begin{align*}
  &\left \| \eps \mapsto {\eps^{m-(n+s)}} |u_{n,\eps}|_{m,p,B} \right \|_{L^q({\R_{>0}},d\eps/\eps)}\\
    &\quad \lesssim \left \| \eps \mapsto {\eps^{m-(n+s)}} 
    \left(\| \DB^{m} u_{n,\eps}\|_p+ [\DB^{m-1} u_{n,\eps}]_{Y,1,p}\right) \right \|_{L^q({\R_{>0}},d\eps/\eps)}
    \intertext{(by Minkowski inequality)}
    &\quad \lesssim \sum_{i = 0}^n{\left \| \eps \mapsto {\eps^{m-(n+s)}} \left(\|I^{(n-i,{m-i})}_{\eps} \DB^i u\|_p +[I^{(n-i,{m-i-1})}_{\eps} \DB^i u]_{Y,1,p}\right) \right\|_{L^q({\R_{>0}},d\eps/{\eps})}}
    \intertext{(by \eqref{eq:estim1_b}-\eqref{eq:estim2_b})}
&\quad \lesssim \sum_{i = 0}^n {|\DB^i u|_{n-i+s,p,q}} {\lesssim | u |_{n+s,p,q}}.
\end{align*}

\medskip 
\noindent \textit{PART 3: Proof of the preliminary estimates \eqref{eq:estim1_b}-\eqref{eq:estim2_b}}.
We only show \eqref{eq:estim2_b} as \eqref{eq:estim1_b} can be proved by similar arguments and it is easier. 
First note that, for any $l\in \N_0$, $\eps>0$ we have
\begin{equation}
 \DB^l\phi(D_{\eps^{-1}}z)=\eps^{-l}( \DB^l\phi)(D_{\eps^{-1}}z), \quad z\in \R^{N+1}
\end{equation}
by \eqref{homogeneity}. Set ${l=m-i-1}$, then we have
\begin{align*}
    \left \| \eps \mapsto {\eps^{m-(n+s)}}[I^{(n-i,l)}_{\eps} \DB^i u]_{Y,1,p} \right\|_{L^q({\R_{>0}},d\eps/{\eps})}^q 
    \lesssim \int_{0}^{{\infty}}\eps^{({i+1-n-s})q} (S_1 + S_2)\frac{d\eps}{\eps},
\end{align*}
where
\begin{align*}
  S_1 := &  {\left( \int_{\R} \int_{\R^{N+1}} 
    \bigg| \int_{\R^{N+1}}\left(T_{n-i} \DB^i u(\z,z) - \DB^i u(z) \right)
    \Delta\phi^l(z,\z,h) \frac{d\z}{\eps^{\Q}} \bigg|^p dz \frac{dh}{|h|^{1+\frac{p}{2}}}\right)^{\frac{q}{p}}}\\ 
    S_2 := &{\left(\int_{\R}  \int_{\R^{N+1}} 
    \bigg| \int_{\R^{N+1}}  \Delta T^{n,i}u(z,\z,h)
    (\DB^{l} \phi ) (D_{\eps^{-1}}(\z^{-1} \circ z)) \frac{d\z}{\eps^{\Q}} \bigg |^p dz 
    \frac{dh}{|h|^{1+\frac{1}{2}q}}\right)^{\frac{q}{p}}} ,
\end{align*}
with
\begin{align*}
\Delta\phi^l_{\eps^{-1}}(z,\z,h) &:= (\DB^{{l}} \phi )(D_{\eps^{-1}}(\z^{-1} \circ e^{hY}z)) - ( \DB^{l} \phi ) (D_{\eps^{-1}}(\z^{-1} \circ z)),\\
\Delta T^{n,i}u(z,\z,h)&:= T_{n-i} \DB^i u(\z,e^{hY}z) - \DB^i u(e^{hY}z) - \left(T_{n-i} \DB^i u(\z,z) - \DB^i u(z) \right).
\end{align*}

By Minkowski's integral inequality
\begin{align*}
   S_{1}\leq & 
    \left(\int_{\R} \bigg( \int_{\R^{N+1}} \bigg( \int_{\R^{N+1}} \left| T_{n-i} \DB^i u(\z,z) - \DB^i u(z) \right|^p  
    |\Delta\phi^l_{\eps^{-1}}(z,\z,h)|^p dz \bigg)^{\frac{1}{p}} \frac{d\z}{\eps^{\Q}} \bigg )^p \frac{dh}{|h|^{1+\frac{p}{2}}}\right)^{\frac{q}{p}} \\
    \intertext{(by a change of variables $\eta=\z^{-1}\circ z$)}
    \leq &\left(\int_{\R}\bigg( \int_{\R^{N+1}} \| T_{n-i} \DB^i u(\cdot \circ \eta^{-1},\cdot) - \DB^i u \|_p  
    \big|\Delta\phi^l_{\eps^{-1}} (\eta,0,h)\big| \frac{d\eta}{\eps^{\Q}} \bigg)^p \frac{dh}{|h|^{1+\frac{p}{2}}}\right)^{\frac{q}{p}} \\
    \intertext{(by a further application of Minkowski's inequality)}
    \leq &\left(\int_{\R^{N+1}} \left(\int_{\R}
    \| T_{n-i} \DB^i u(\cdot \circ \eta^{-1},\cdot) - \DB^i u \|_p^p 
    \big|\Delta\phi^l_{\eps^{-1}} (\eta,0,h)\big|^p\frac{dh}{|h|^{1+\frac{p}{2}}}\right)^{\frac{1}{p}} \frac{d\eta}{\eps^{\Q}}\right)^q\\
    =&\left(\int_{\R^{N+1}} \| T_{n-i} \DB^i u(\cdot \circ \eta^{-1},\cdot) - \DB^i u \|_p \left(\int_{\R}
    \big|\Delta\phi^l_{\eps^{-1}} (\eta,0,h)\big|^p\frac{dh}{|h|^{1+\frac{p}{2}}}\right)^{\frac{1}{p}} \frac{d\eta}{\eps^{\Q}}\right)^q
\end{align*}
Observe that on the set $A=\{\|\eta\|_B<\|e^{hY}\eta\|_B\}$, the support of $\Delta\phi^l (\eta,0,h)$ is included in 
$\|D_{\eps^{-1}}\eta\|_B=\eps^{-1}\|\eta\|_B\le 1$ and on $A^c$ it is included in $\|D_{\eps^{-1}}\eta\|_B\le 1$. 
Therein we have, by Lemma \ref{lem0},
$$\|D_{\eps^{-1}}\eta\|_B=\|e^{hY}D_{\eps^{-1}}\eta\circ (-h,0) \|_B\le m(\|D_{\eps^{-1}}\eta\|_B+|h|^{\frac 12})\le 2m,$$
therefore the support of $\Delta\phi^l (\eta,0,h)$ is always included in $\|\eta\|_{B}\le 2m\eps$.
Moreover, by \eqref{homogeneity} we have, 
\begin{align*}
&\int_{\R}\big|\Delta\phi^l_{\eps^{-1}} \big(D_{\eps^{-1}}\eta,0,\frac{h}{\eps^2}\big)\big|^p\frac{dh}{|h|^{1+\frac{p}{2}}}\\
&\qquad=\left(\int_{|h| \leq \eps^2}+\int_{|h|\ge \eps^2}\right)\big|\Delta\phi^l_{\eps^{-1}} (\eta,0,h)\big|^p\frac{dh}{|h|^{1+\frac{p}{2}}}
\intertext{(denoting with $\chi_{\{\|\eta\|_B\le 2m\eps\}}$  the indicator function of the set $\|\eta\|_B\le 2m\eps$)}
&\qquad \lesssim \chi_{\{\|\eta\|_B\le 2m\eps\}}\left(\|Y\DB^l\phi\|^p_{\infty}\int_{|h| \leq \eps^2}\frac{dh}{\eps^{2p}|h|^{1-\frac{p}{2}}}+
\|\DB^l\phi\|^p_{\infty}\int_{|h| \geq \eps^2}\frac{dh}{|h|^{1+\frac{p}{2}}}\right)\\
&\qquad \lesssim \chi_{\{\|\eta\|_B\le 2m\eps\}} \eps^{-p}.
\end{align*}
{Therefore, by H\"older inequality we get
\begin{align*}
&\int_{0}^{\infty}\eps^{{(i+1-n-s)}q} S_{1}\frac{d\eps}{\eps}\\
&\qquad\lesssim \int_{0}^{\infty} \left ( \int_{\|\eta\|_B \le 2m \eps} 
\| T_{n-i} \DB^i u(\cdot \circ \eta^{-1},\cdot) - \DB^i u \|_p
\frac{d\eta}{\eps^{\Q}} \right )^q \frac{d\eps}{\eps^{1+(n+s-i)q}} \\
&\qquad\leq \int_{0}^{\infty} \int_{\|\eta\|_B \le 2m \eps}
\| T_{n-i} \DB^i u(\cdot \circ \eta^{-1},\cdot) - \DB^i u \|_p^q \frac{d\eta}{\eps^{\Q}} 
\left ( \int_{\|\eta\|_B \le 2m \eps} \frac{d\eta}{\eps^{\Q}} \right )^{\frac{q}{q'}}
 \frac{d\eps}{\eps^{1+(n+s-i)q}} \\
&\qquad\lesssim \int_{\frac{1}{2m}\|\eta\|_B}^{\infty} \int_{\R^{N+1}}
\| T_{n-i} \DB^i u(\cdot \circ \eta^{-1},\cdot) - \DB^i u \|_p^q d\eta \frac{d\eps}{\eps^{1+\Q+(n+s-i)q}} \\
& \qquad \lesssim \int_{\R^{N+1}} \| T_{n-i} \DB^i u(\cdot \circ \eta^{-1},\cdot) - \DB^i u \|_p^q 
\frac{d\eta}{\|\eta\|_{B}^{\Q+(n+s-i)q}} \lesssim {| \DB^i u|_{n-i+s,p,q}^q},
\end{align*}
}
where the last inequality stems from Proposition \ref{prop:prop1} and Theorem \ref{th:Taylor}. 

Next we consider $S_2$. With similar computation as before we have
\begin{align*}
   S_{2}\leq & 
    \left(\int_{\R} \bigg( \int_{\R^{N+1}} \bigg( \int_{\R^{N+1}} \left| \Delta T^{n,i}u(z,\z,h) \right|^p  
    | (\DB^{l} \phi ) (D_{\eps^{-1}}(\z^{-1} \circ z))|^p dz \bigg)^{\frac{1}{p}} \frac{d\z}{\eps^{\Q}} \bigg )^p \frac{dh}{|h|^{1+\frac{p}{2}}}\right)^{\frac{q}{p}} \\
      \leq &\left(\int_{\R}\bigg( \int_{\|\eta\|_B\le \eps} \| \Delta T^{n,i}u(\cdot,\cdot\circ \eta^{-1},h) \|_p  
    \big|(\DB^{l} \phi ) (D_{\eps^{-1}}\eta)\big| \frac{d\eta}{\eps^{\Q}} \bigg)^p \frac{dh}{|h|^{1+\frac{p}{2}}}\right)^{\frac{q}{p}}
    \\ \lesssim & S_{21}^q+S_{22}^q,
\end{align*}
where 
\begin{align*}
S_{21}:&=\left(\int_{|h| \geq \eps^2}\bigg( \int_{\|\eta\|_B\le \eps} \| \Delta T^{n,i}u(\cdot,\cdot\circ \eta^{-1},h) \|_p  
    \big|(\DB^{l} \phi ) (D_{\eps^{-1}}\eta)\big| \frac{d\eta}{\eps^{\Q}} \bigg)^p \frac{dh}{|h|^{1+\frac{p}{2}}}\right)^{\frac{1}{p}}
\intertext{(by Minkowski's inequality)}  
&\lesssim     \int_{\|\eta\|_B\le \eps} \left(\int_{|h| \geq \eps^2}
  \| \Delta T^{n,i}u(\cdot,\cdot\circ \eta^{-1},h) \|^p_p\frac{dh}{|h|^{1+\frac{p}{2}}} \right)^{\frac{1}{p}}
   \big|(\DB^{l} \phi ) (D_{\eps^{-1}}\eta)\big| \frac{d\eta}{\eps^{\Q}}\\
 &\lesssim     \int_{\|\eta\|_B\le \eps} \left(\int_{|h| \geq \eps^2}
 \| T_{n-i}\DB^i u(\cdot\circ \eta^{-1},\cdot)- \DB^i u \|^p_p\frac{dh}{|h|^{1+\frac{p}{2}}} \right)^{\frac{1}{p}}
   \big|(\DB^{l} \phi ) (D_{\eps^{-1}}\eta)\big| \frac{d\eta}{\eps^{\Q}}
   \intertext{(evaluating the $dh$ integral)}  
  &\lesssim  \int_{\|\eta\|_B\le \eps} \eps^{-1}\| T_{n-i}\DB^i u(\cdot\circ \eta^{-1},\cdot)- \DB^i u \|_p
   \big|(\DB^{l} \phi ) (D_{\eps^{-1}}\eta)\big| \frac{d\eta}{\eps^{\Q}},
   \intertext{and}
  S_{22}:&=\left(\int_{|h| \leq \eps^2}\bigg( \int_{\|\eta\|_B\le \eps} \| \Delta T^{n,i}u(\cdot,\cdot\circ \eta^{-1},h) \|_p  
    \big|(\DB^{l} \phi ) (D_{\eps^{-1}}\eta)\big| \frac{d\eta}{\eps^{\Q}} \bigg)^p \frac{dh}{|h|^{1+\frac{p}{2}}}\right)^{\frac{1}{p}}.
\end{align*}
In particular we have
\begin{align*}
&\int_{0}^{\infty}\eps^{{(i+1-n-s)}q} S_{21}^q\frac{d{\eps}}{{\eps}}\\
 &\qquad \lesssim \int_{0}^{\infty}\eps^{{(i-n-s)}q}\left(\int_{\|\eta\|_B\le {\eps}} \| T_{n-i}\DB^i u(\cdot\circ \eta^{-1},\cdot) - \DB^i u \|_p \big|(\DB^{l} \phi ) (D_{\eps^{-1}}\eta)\big| \frac{d\eta}{\eps^{\Q}}\right)^q \frac{d{\eps}}{{\eps}}
 \intertext{(by {H\"older}'s inequality)}
 &\qquad \lesssim \int_{ \|\eta\|_{B}}^{\infty} \int_{\R^{N+1}} \| T_{n-i} \DB^i u(\cdot \circ \eta^{-1},\cdot) - \DB^i u \|_p^q d\eta \frac{d{\eps}}{\eps^{1+\Q+(n+s-i)q}} \\
& \qquad  \lesssim \int_{\R^{N+1}} \| T_{n-i} \DB^i u(\cdot \circ \eta^{-1},\cdot) - \DB^i u \|_p^q 
\frac{d\eta}{\|\eta\|_{B}^{{\Q+(s+n-i)q}}} \lesssim {| \DB^i u|_{n+s-i,p,q}^q}.
\end{align*}
It remains to check the term in {$S_{22}$}. We only consider the case $n-i\ge 2$ the other ones being simpler: indeed
a simple computation shows that $T_0\DB^iu(\z,z)=u(\z)=T_0(\z,e^{hY}z)$ and $T_1\DB^iu(\z,z)=T_1(\z,e^{hY}z)$ for any $z,\z\in\R^{N+1}$, $h\in \R$, which yields
$$\Delta T^{n,i}u(z,z\circ \eta{-1},h)=\DB^iu(z)-\DB^iu(e^{hY}z), \quad  n=i, \; n=i+1.$$
On the other hand, by the Mean Value Theorem along the vector field $Y$, Proposition \ref{prop:prop1} and Lemma \ref{LEM}, we have for $n-i\ge 2$
\begin{align*}
\| \Delta T^{n,i}u(\cdot,\cdot\circ \eta^{-1},h) \|_p\lesssim 
|h| \left\| 
T_{n-i-2}Y\DB^i u(\cdot\circ \eta^{-1},e^{\bar{h} Y}\cdot)- Y\DB^i u(e^{\bar{h} Y}\cdot)\right\|^p_p
\end{align*}
for some $\bar{h}\in [0,h]$. Notice that
  $$z\circ \eta^{-1}=e^{\bar{h}Y}z\circ (-\bar{h},0)\circ  \eta^{-1}=e^{\bar{h}Y}z\circ (e^{\bar{h}Y} \eta)^{-1}.$$
Then, after the change of variables $\bar z=e^{\bar{h}Y}z$ inside the $L^p$ norm
\begin{align*}
\| \Delta T^{n,i}u(\cdot,\cdot\circ \eta^{-1},h) \|_p\lesssim |h|  \left\|
T_{n-i-2}Y\DB^i u(\cdot\circ (e^{\bar{h} Y}\eta)^{-1},\cdot)- Y\DB^i u \right\|_p.
\end{align*}
Therefore we have
\begin{align*}
S_{22}&\lesssim \left(\int_{|h| \leq {\eps}^2}\bigg( \int_{\|\eta\|_B\le \eps} 
\left\|T_{n-i-2}Y\DB^i u(\cdot\circ (e^{\bar{h} Y}\eta)^{-1},\cdot)- Y\DB^i u \right\|_p  
    \big|(\DB^{l} \phi ) (D_{\eps^{-1}}\eta)\big| \frac{d\eta}{\eps^{\Q}} \bigg)^p \frac{dh}{|h|^{1-\frac{p}{2}}}\right)^{\frac{1}{p}}
    \intertext{(by the change of variable $\bar{\eta}=e^{\bar{h} Y}\eta$)}
&\lesssim \left(
 \int_{|h| \leq \eps^2}
\left( \int_{\|e^{-\bar{h} Y}\bar{\eta}\|_B\le {\eps}} 
\left\|T_{n-i-2}Y\DB^i u(\cdot\circ \bar{\eta}^{-1},\cdot)- Y\DB^i u \right\|_p  
\big|(\DB^{l} \phi ) (D_{\eps^{-1}}e^{-\bar{h} Y}\bar{\eta})\big|
\frac{d\bar\eta}{\eps^{\Q}} \right)^p \frac{dh}{|h|^{1-\frac{p}{2}}} \right)^{\frac{1}{p}}   
\end{align*}
Note that on the current integration set we have, by {Lemma \ref{lem0}}
$$\|\bar\eta\|_B=\|e^{-\bar{h} Y}\bar{\eta}\circ (\bar h,0)\|_B
\le m(|\bar h|^{\frac 12}+\|e^{-\bar{h} Y}\bar{\eta}\|_B)\le 2m\eps. $$
Then, again by Minkowski inequality we can write 
\begin{align*}
S_{22}^q &\lesssim \left( \int_{\|{\eta}\|_B\le 2m \eps} \left( \int_{|h| \leq \eps^2} 
\left\|T_{n-i-2}Y\DB^i u(\cdot\circ {\eta}^{-1},\cdot)- Y\DB^i u \right\|^p_p  
\big|(\DB^{l} \phi ) (D_{\eps^{-1}}e^{-\bar{h} Y}{\eta})\big|^p
\frac{dh}{|h|^{1-\frac{p}{2}}} \right)^{\frac{1}{p}} \frac{d\eta}{\eps^{\Q}} \right)^q\\
\intertext{(since $\DB^{l} \phi$ is bounded and evaluating the $dh$ integral)}
&\lesssim \eps^q \left( \int_{\|{\eta}\|_B\le 2m \eps} \left\|T_{n-i-2}Y\DB^i u(\cdot\circ {\eta}^{-1},\cdot)- Y\DB^i u \right\|_p \frac{d\eta}{\eps^{\Q}} \right)^q\\
\end{align*}
Finally, reasoning as for the term $S_1${
\begin{align*}
&\int_{0}^{\infty}\eps^{{(i+1-n-s)}q} S_{22}^q\frac{d\eps}{\eps}\\
&\qquad \lesssim \int_{0}^{\infty} \left ( \int_{\|\eta\|_B \le 2m\eps} 
\| T_{n-i-2} Y\DB^i u(\cdot \circ \eta^{-1},\cdot) - Y\DB^i u \|_p  \frac{d\eta}{\eps^{\Q}} \right )^q 
\frac{d\eps}{\eps^{1+{(n+s-i-2)}q}} \\
\intertext{(by H\"older's inequality)}
&\qquad \lesssim \int_{\frac{1}{2m} \|\eta\|_B}^{\infty} \int_{\R^{N+1}} \| T_{n-i-2} Y\DB^i u(\cdot \circ \eta^{-1},\cdot) - Y\DB^i u \|_p^q d\eta \frac{d\eps}{\eps^{1+\Q+{(n+s-i-2)}q}} \\
& \qquad  \lesssim \int_{\R^{N+1}}\| T_{n-i-2} Y\DB^i u(\cdot \circ \eta^{-1},\cdot) - Y\DB^i u \|_p^q \frac{d\eta}{\|\eta\|_{B}^{\Q+(n+s-i-2)q}} \lesssim {| Y\DB^i u|_{n+s-i-2,p,q}^q}.
\end{align*}
}
The proof is complete.
\endproof

\subsection{Proof of Theorem \ref{th:interpolation}}
In this Section we denote $\tilde{\La}^{n+s}_{p,q,B}:=\big(L^p,W^{n+1,p}_{B}\big)_{\frac{n+s}{n+1},q}$  and prove that 
$\tilde{\La}^{n+s}_{p,q,B}={\La}^{n+s}_{p,q,B}$.


\subsubsection*{Embedding $\tilde{\La}^{n+s}_{p,q,B} {\subseteq} {\La}^{n+s}_{p,q,B}$}
We proceed by a two step induction on $n$.

\medskip

\noindent
\textit{Case} $n=0$: for $u \in \tilde{\La}^{s}_{p,q,B}$, there exist $u_{1,\eps}\in L^p $ and $u_{2,\eps}\in W^{1,p}_B$, 
$\eps\ge 0$, such that 
$u=u_{1,\eps}+u_{2,\eps}$ and
$$\|u_{1,\eps}\|_p+\eps \|u_{2,\eps}\|_{W^{1,p}_B}\le 2K(\eps,u,L^p,W^{1,p}_B),$$ where $\eps\mapsto \eps^{-s}K(\eps,u,L^p,W^{1,p}_B)\in L^q(\R_{>0},{d\eps}/\eps)$ by definition. By Theorem 4.1 in \cite{MR4700191}) we have 
\begin{align*}
\|u(e^{h Y}\cdot)-u\|_p&\le \|u_{1,|h|^{\frac 12}}(e^{h Y}\cdot)-u_{{1,|h|^{\frac 12}}} \|_p+\|u_{2,|h|^{\frac 12}}(e^{h Y}\cdot)-u_{{2,|h|^{\frac 12}}} \|_p\\
&\le 2\|u_{1,|h|^{\frac 12}}\|_p+|h|^{\frac 12}[u_{2,|h|^{\frac 12}}]_{Y,1,p,p}\le 4 K(|h|^{\frac 12},u,L^p,W^{1,p}_B).
\end{align*}
Then we have
\begin{align*}
    [u]^q_{Y,s,p,q} &=\int_{0}^{\infty} \|u(e^{h Y}\cdot)-u\|_p^q \frac{dh}{|h|^{1+\frac{s}{2}q}}+\int_{-\infty}^{0} \|u(e^{h Y}\cdot)-u\|_p^q \frac{dh}{|h|^{1+\frac{s}{2}q}}
    \intertext{(by the change of variable $z'=e^{hY}z$ in the inner norm of the second term and recalling Remark \ref{r1})}
    &=\int_{0}^{\infty} \|u(e^{h Y}\cdot)-u\|_p^q \frac{dh}{|h|^{1+\frac{s}{2}q}}+\int_{-\infty}^{0} \|u-u(e^{-h Y}\cdot)\|_p^q \frac{dh}{|h|^{1+\frac{s}{2}q}}\\
    &\lesssim \int_{0}^{\infty} \|u(e^{h Y}\cdot)-u\|_p^q \frac{dh}{|h|^{1+\frac{s}{2}q}}\\
    &\lesssim \left\|h\mapsto h^{-\frac{s}{2}}K(h^{\frac{1}{2}},u,L^p,W^{1,p}_B) \right\|^q_{L^{q}(\R_{>0},\frac{dh}{h})}
    =\|u\|^q_{\tilde{\La}^{s}_{p,q,B}}.
\end{align*}
Similarly, for any $i=1,\dots,d$ we have 
\begin{align*}
    \|u(e^{h \p_{x_i}}\cdot)-u\|_p
    &\le 2\|u_{1,|h|}\|_p+|h|\|\p_{x_i}u_{2,|h|}\|_p \le 4 K(|h|,u,L^p,W^{1,p}_B).
\end{align*}
and therefore 
\begin{align*}
    [u]_{\p_{x_i},s,p,q} &\lesssim
    \left\| h\mapsto h^{-{s}}K(h,u,L^p,W^{1,p}_B) \right\|_{L^{q}(\R_{>0},\frac{dh}{h})}
    =\|u\|_{\tilde{\La}^{s}_{p,q,B}}.
\end{align*}
On the other hand, clearly $\tilde{\La}^{s}_{p,q,B}\subseteq L^p$ so that $\|u\|_p\le \|u\|_{\tilde{\La}^{s}_{p,q,B}}$. This concludes the proof for $n=0$.

\medskip
\noindent \textit{Case} $n=1$: first notice that, for any $n\in \N_0$, and $X\in \{Y,\p_{x_1},\dots, \p_{x_d}\}$, we have the trivial estimate 
\begin{equation}\label{eq:interp_e1}
K(\l,Xu,W^{n,p}_B,W^{n+1,p}_B)\le K(\l,u,W^{n+m_X,p}_B,W^{n+1+m_X,p}_B), 
\end{equation}
for any $u\in W^{n+m_X,p}_B+W^{n+1+m_X,p}_B$. 
Next, we need to prove that, for any $u \in \tilde{\La}^{1+s}_{p,q,B}$ we have 
\begin{align*}
[u]_{Y,{1+s},p,q} \lesssim \|u\|_{\tilde{\La}^{1+s}_{p,q,B}}, \quad 
{|\p_{x_i}u|_{s,p,q}\lesssim \|u\|_{\tilde{\La}^{1+s}_{p,q,B}}}.
\end{align*}
To prove the first inclusion we use the characterization 
$\tilde{\La}^{1+s}_{p,q,B}=(L^{p},W^{2,p}_B)_{\frac{1+s}{2},q}$. Taking 
$u_{1,\eps}\in L^p $ and $u_{2,\eps}\in W^{2,p}_B$, $\eps\ge 0$, such that 
$u=u_{1,\eps}+u_{2,\eps}$ and
$$\|u_{1,\eps}\|_p+\l \|u_{2,\eps}\|_{W^{2,p}_B}\le 2K(\eps,u,L^p,W^{2,p}_B),$$ with $\eps\mapsto \eps^{-\frac{1+s}{2}}K(\eps,u,L^p,W^{2,p}_B)\in L^q(\R_{>0},{d\eps}/\eps)$, we have, analogously to step 1, 
\begin{align*}
\|u(e^{h Y}\cdot)-u\|_p
&\le 2\|u_{1,|h|}\|_p+|h|\|Yu_{2,|h|}\|_p \le 4 K(|h|,u,L^p,W^{2,p}_B),
\end{align*}
and therefore  
\begin{align*}
[u]_{Y,{1+s},p,q} \lesssim
 \left\|h\mapsto h^{-\frac{1+s}{2}}K(h^{\frac{1}{2}},u,L^p,W^{1,p}_B) \right\|_{L^{q}(\R_{>0},\frac{dh}{h})}
 = {\|u\|_{\tilde{\La}^{1+s}_{p,q,B}}}.
\end{align*}
To prove the second inequality it is more convenient to use the characterization by reiteration 
$\tilde{\La}^{1+s}_{p,q,B}=(W^{1,p}_B,W^{2,p}_B)_{s,q}$: indeed, by the previous step and \eqref{eq:interp_e1}, for any $i=1,\cdots d$ we have 
\begin{align*}
|\p_{x_i}u|_{s,p,q}
&\lesssim \|\p_{x_i}u\|_{\tilde{\La}^{s}_{p,q,B}}\\
&\lesssim\left\|h\mapsto h^{-\frac{s}{2}}K(h^{\frac{1}{2}},\p_{x_i}u,L^p,W^{1,p}_B) \right\|_{L^{q}(\R_{>0},\frac{dh}{h})}\\
&{\lesssim\left\|h\mapsto h^{-\frac{s}{2}}K(h^{\frac{1}{2}},u,W^{1,p}_B,W^{2,p}_B) \right\|_{L^{q}(\R_{>0},\frac{dh}{h})}
\lesssim\|u\|_{\tilde{\La}^{1+s}_{p,q}(B)}},
\end{align*}
where in the last equality we used the change of variable $h'=h^{\frac{1}{2}}$.

\medskip
\noindent \textit{Case} $n\ge 2$: Assume now the general embedding holds for some $n\in \N$ and $n+1$, and let $u\in \La^{n+2+s}_{p,q,B}$. Then by \eqref{eq:interp_e1} we have
\begin{align*}
\|u\|_{\La^{n+2+s}_{p,q,B}} &= \|u\|_p + |Yu|_{n+s,p,q}
+ \sum_{i=1}^{d}|\p_{x_i}u|_{n+1+s,p,q}\\
& \lesssim \|u\|_p+\|Yu\|_{\tilde{\La}^{n+s}_{p,q,B}}+\sum_{i=1}^{d}\|\p_{x_i}u\|_{\tilde{\La}^{n+1+s}_{p,q,B}}
\intertext{(characterizing by reiteration $\tilde{\La}^{n+s}_{p,q,B}=(W^{n,p}_B,W^{n+1,p}_B)_{s,q}$ and $\tilde{\La}^{n+1+s}_{p,q,B}=(W^{n+1,p}_B,W^{n+2,p}_B)_{s,q}$)}
&\lesssim \|u\|_{p}+
\left\|h\mapsto h^{-s}K(h^s,Yu,W^{n,p}_B,W^{n+1,p}_B) \right\|_{L^{q}(\R_{>0},\frac{dh}{h})}\\
&\qquad + \sum_{i=1}^d 
\left\|h\mapsto h^{-s}K(h^s,\p_{x_i}u,W^{n+1,p}_B,W^{n+2,p}_B) \right\|_{L^{q}(\R_{>0},\frac{dh}{h})}
\intertext{(by \eqref{eq:interp_e1} and again by reiteration)}
&\lesssim \|u\|_{(W^{n+2,p}_B,W^{n+3,p}_B)_{s,q}}\lesssim \|u\|_{\tilde{\La}^{n+2+s}_{p,q,B}}.
\end{align*}
The proof is complete.
\endproof


\subsubsection*{Embedding ${\La}^{n+s}_{p,q,B} \subseteq \tilde{\La}^{n+s}_{p,q,B}$}

Embedding ${\La}^{n+s}_{p,q,B} \subseteq \tilde{\La}^{n+s}_{p,q,B}$ is a direct consequence of Proposition \ref{prop:approx}.
Indeed, for any $\eps>0$ we have 
$$K(h,u;L^p,W^{n+1,p}_B)\le \|u-u_{n,\eps}\|_p+\l\|u_{n,\eps}\|_{W^{n+1,p}_B}$$
Then,  by \eqref{eq:approx1} and \eqref{eq:approx2},  for any $u\in \La^{n+s}_{p,q,B}$ and taking $\eps=h^{\frac{1}{n+1}}$ we get
\begin{align*}
\|u\|_{\tilde{\La}^{n+s}_{p,q,B}}&=\left\|h\mapsto h^{-\frac{n+s}{n+1}}K(h,u;L^p,W^{n+1,p}_B) \right\|_{L^q(\R_{>0},\frac{dh}{h})}\\
&{\leq}\left\|h\mapsto h^{-\frac{n+s}{n+1}}\left(\Big\|u-u_{n,h^{\frac{1}{n+1}}}\Big\|_p + h \Big\|u_{n,h^{\frac{1}{n+1}}}\Big\|_{W^{n+1,p}_B}\right)\right\|_{L^q(\R_{>0},\frac{dh}{h})}
\intertext{(by the change of variable $h'=h^{\frac{1}{n+1}}$)}
&\lesssim \left\|h\mapsto h^{-{(n+s)}}\|u-u_{n,h}\|_p  
+h^{1-s}\|u_{n,h}\|_{W^{n+1,p}_B}\right\|_{L^q(\R_{>0},\frac{dh}{h})}\lesssim \|u\|_{{\La}^{n+s}_{p,q,B}},
\end{align*}
where we exploited Proposition \ref{prop:approx} in the last inequality.
The proof is complete.
\endproof

\begin{corollary}[Interpolation inequalities]
Let$s,s_1\in (0,1)$, $n,n_1\in \N_0$ with $n_1>n$ and $p,q,q_1\ge 1$. Then there  exists a positive constant $c=c(s,s_1,n,n_1,p,q,q_1,B)$ such that, for any $u\in \Lambda^{n_1+s_1}_{p,q_1,B}$ and $\eps>0$, we have 
\begin{equation}\label{eq:cor2}
|u|_{n+s,p,q}\le c \left(\eps^{-\frac{-n-s}{n_1+s_1-n-s}}\|u\|_p +\eps |u|_{n_1+s_1,p,q_1}\right).
\end{equation}
In particular, optimizing w.r.t $\eps$ we obtain 
\begin{equation}\label{eq:cor1}
|u|_{n+s,p,q}\le c \|u\|_p^{1-\frac{n+s}{n_1+s_1}} \left(|u|_{n_1+s_1,p,q_1}\right)^{\frac{n+s}{n_1+s_1}}.
\end{equation}
\end{corollary}
\proof By Theorem \ref{th:interpolation} and reiteration, the space $\Lambda^{n+s}_{p,q,B}$ is intermediate between $L^p$ and $\Lambda^{n_1+s_1}_{p,q_1,B}$. Thus we get \eqref{eq:cor2} by scaling (cf. Lemma \ref{scaling}). 
\endproof

\section{Embeddings of $\La^{n+s}_{p,q,B}$ spaces}\label{sec:emb}

\subsection{Proof of Theorem \ref{th:embeddings}}

The proof of Theorem \ref{th:embeddings} fully exploits the interpolation characterization of the Besov spaces of Theorem \ref{th:interpolation}, as well as the embeddings for the intrinsic Sobolev spaces $W^{m,p}_{B}$ proved in \cite{MR4700191} and the interpolation characterization of intrinsic H\"older spaces proved in \cite{MR4694350}. For completeness, we include their statements below. 

\begin{theorem}[Real interpolation of intrinsic H\"older spaces]\label{th:Holder_interp}
Let $n_1,n_2\in\N_0$ and $\a_1,\a_2\in [0,1]$, with $n_1+\a_1\le n_2+\a_2$, and let $\theta\in (0,1)$. Let $n\in\N_0$, $n_1\le n\le n_2$ such that
\begin{equation}\label{eTHbis}
\frac{n-(n_1+\a_1)}{(n_2+\a_2)-(n_1+\a_1)}<\theta <\frac{n+1-(n_1+\a_1)}{(n_2+\a_2)-(n_1+\a_1)}.
\end{equation}
If $\a:=(n_1+\a_1)+\theta[n_2+\a_2-(n_1+\a_1)]-n \notin \{0,1\}$
then it holds that
\begin{equation}\label{eTH}
\left(C^{n_1,\a_1}_B,C^{n_2,\a_2}_B\right)_{\theta,\infty}=C^{n,\alpha}_B
\end{equation}
\end{theorem}

\begin{theorem}[Embeddings for natural order intrinsic Sobolev spaces]\label{th:Sob_emb}
Let $k\in \N$, $n\in \N_0$ and ${p\in (1,\infty)}$.
\begin{itemize}
  \item[i)] If $kp<\Q$, then
\begin{equation}\label{Sobolev_n}
  W^{n+k,p}_B\subseteq W_B^{n,q}, \qquad p\le q\le p^{\ast}_{k},
  \qquad \frac{1}{p^{\ast}_{k}}=\frac{1}{p}-\frac{k}{\Q}.
\end{equation}
In particular, $W^{k,p}_B\subseteq L^{q}$ for $p\le q\le p^{\ast}_{k}$;
\item[ii)] if $np>\Q$, then
\begin{equation}\label{Morrey_n}
  W^{n,p}_B\subseteq C_B^{\left[n-\frac{\Q}{p}\right],\left\{n-\frac{\Q}{p}\right\}}
\end{equation}
where $[\cdot]$ and $\{\cdot\}$ denote the integer and fractional part respectively.
\end{itemize}
\end{theorem}

\noindent We need one last lemma.
\begin{lemma}\label{lem:emb}For $m\in \N$ and $p>\frac{\bf d}{m}$ we have 
\begin{equation}\label{eq:emb_0}
\left(L^{p}, W^{m,p}_{B}\right)_{\theta,1}\subseteq L^{\infty}, \qquad \theta=\frac{\bf d}{m p}.
\end{equation}
\end{lemma}
\proof 
By \cite{MR4700191}, Theorems 1.1 (ii) and Theorem 7.1 (ii), it is already known that $W^{m,p}_{ B}\subseteq L^{\infty}$. Therefore, using the characterization of Definition 2.7, we have, for any $u\in  W^{m,p}_{ B}$: 
\begin{equation}\label{eq:emb_1}
\|u\|_{\infty}\lesssim \|u\|_p+|u|_{m,p, B}.
\end{equation}
Next we exploit a scaling argument: applying \eqref{eq:emb_1} to $u(D_\l)$ and exploiting Lemma 2.12 (in \cite{MR4700191}) we get
\begin{align*}
\|u\|_{\infty}&\lesssim \|u(D_\l)\|_p+|u(D_\l)|_{m,p, B}
\lesssim \l^{-\frac{\bf d}{p}} \|u\|_p+\l^{m-\frac{\bf d}{p}}|u|_{m,p, B}.
\end{align*}
Optimizing w.r.t. $\l$ we directly infer
\begin{equation}\label{eq:emb_2}
\|u\|_{\infty}\lesssim \|u\|_p^{1-\theta}|u|_{m,p, B}^{\theta}\lesssim \|u\|_p^{1-\theta}\|u\|_{W^{m,p}_{B}}^{\theta},\qquad \theta=\frac{\bf d}{m p}.
\end{equation}
The equivalence between \eqref{eq:emb_0} and \eqref{eq:emb_2} is standard, see for instance \cite{MR3753604}, Proposition 1.20.
\endproof

\noindent We are ready to prove Theorem \ref{th:embeddings}.

\proof[Proof of Theorem \ref{th:embeddings}]

\textit{i)} For $(n+s)p={\bf d}$, taking $m\in \N$, $m>\frac{\bf d}{p}$ we have
$$\Lambda^{n+s}_{p,1,B}=\left(L^{p}, W^{m,p}_{B}\right)_{\frac{n+s}{m},1}{=} \left(L^{p}, W^{m,p}_{B}\right)_{\frac{\bf d}{m p},1}\subseteq L^{\infty},$$
where the last inclusion holds by Lemma \ref{lem:emb}.


\textit{ii)} Assume now $(n+s)p<{\bf d}$. Let $n_1\in \N$ and $s_1\in [0,1)$ such that $n_1+s_1=\frac{\bf d}{p}$. 
By the reiteration theorem and Lemma \ref{lem:emb} we directly infer
$$\Lambda^{n+s}_{p,q,B}=
\left(L^{p}, \left(L^{p},W^{\frac{\bf d}{p}+1,p}_B\right)_{\frac{\bf d}{{\bf d}+p},1}\right)_{\frac{(n+s) p}{\bf d},q} 
\subseteq \left(L^{p}, L^{\infty}\right)_{\frac{(n+s) p}{\bf d},q}= L^{p',q}, \qquad \frac{1}{p'}=\frac{1}{p}-\frac{n+s}{\bf d}.$$

\textit{iii)} Finally, assume $(n+s)p>{\bf d}$ 
Using again the reiteration theorem we may write, for some  $l\in \N$, $l>n+s$,
\begin{align}\label{eq:emb_3}
\Lambda^{n+s}_{p,q,B}&=\left(\left(L^{p},W^{\frac{\bf d}{p}+1,p}_B\right)_{\frac{\bf d}{{\bf d}+p},1}, W^{l,p}_{B}\right)_{\theta,q}, 
\qquad (1-\theta)\frac{\bf d}{p} +\theta l=n+s.
\intertext{(by Lemma \ref{lem:emb})}
&\subseteq \left(L^{\infty},W^{l,p}_{B}\right)_{\theta,q} \nonumber
\intertext{(by Theorem \ref{th:Sob_emb}, ii)}
&\subseteq \left(L^{\infty},C_{B}^{l-k,k-\frac{\bf d}{p}}\right)_{\theta,q}\subseteq 
\left(L^{\infty},C_{B}^{l-k,k-\frac{\bf d}{p}}\right)_{\theta,\infty},\nonumber
\end{align}
where $k\in \N_0$ is such that $kp >{\bf d} > (k-1)p$.
By Theorem \ref{th:Holder_interp}, we finally infer that
 $\Lambda^{n+s}_{p,q,B}\subseteq C_{B}^{n-k',\a}$,  where $k'$ is the smallest integer such that $(k'+s)p>{\bf d}$ and, recalling \eqref{eq:emb_3}, $$\a=\theta \big(l-\frac{\bf d}{p}\big)-(n-k')=k'+s-\frac{\bf d}{p}.$$
 The proof is complete.
 \endproof
 
 \medskip
 
%
 
%
%


\bibliographystyle{acm}
\bibliography{bib1}
\end{document}